\documentstyle[12pt]{article}

\makeatletter
\def\@sect#1#2#3#4#5#6[#7]#8{\ifnum #2>\c@secnumdepth
    \def\@svsec{}\else
    \refstepcounter{#1}\edef\@svsec{\csname the#1\endcsname.\hskip 1em }\fi
    \@tempskipa #5\relax
    \ifdim \@tempskipa>\z@
    \begingroup #6\relax
    \@hangfrom{\hskip #3\relax\@svsec}{\interlinepenalty \@M #8\par}
    \endgroup
    \csname #1mark\endcsname{#7}\addcontentsline
    {toc}{#1}{\ifnum #2>\c@secnumdepth \else
     \protect\numberline{\csname the#1\endcsname}\fi
           #7}\else
    \def\@svsechd{#6\hskip #3\@svsec #8\csname #1mark\endcsname
          {#7}\addcontentsline
          {toc}{#1}{\ifnum #2>\c@secnumdepth \else
     \protect\numberline{\csname the#1\endcsname}\fi
           #7}}\fi
     \@xsect{#5}}
\def\label#1{\@bsphack\if@filesw {\let\thepage\relax
   \xdef\@gtempa{\write\@auxout{\string
   \newlabel{#1}{{\thesection.\@currentlabel}{\thepage}}}}}\@gtempa
   \if@nobreak \ifvmode\nobreak\fi\fi\fi\@esphack}
\def\@eqnnum{(\thesection.\theequation)}
\def\section{\setcounter{equation}{0} \@startsection {section}{1}{\z@}{-3.5ex
   plus -1ex minus -.2ex}{2.3ex plus .2ex}{\Large\bf}}
\newcount\@minsofar
\newcount\@min
\newcount\@cite@temp
\def\@citex[#1]#2{%
\if@filesw \immediate \write \@auxout {\string \citation {#2}}\fi
\@tempcntb\m@ne \let\@h@ld\relax \def\@citea{}%
\@min\m@ne%
\@cite{%
  \@for \@citeb:=#2\do {\@ifundefined {b@\@citeb}%
    {\@h@ld\@citea\@tempcntb\m@ne{\bf ?}%
    \@warning {Citation `\@citeb ' on page \thepage \space undefined}}%
{\@minsofar\z@ \@for \@scan@cites:=#2\do {%
  \@ifundefined{b@\@scan@cites}%
    {\@cite@temp\m@ne}
    {\@cite@temp\number\csname b@\@scan@cites \endcsname \relax}%
\ifnum\@cite@temp > \@min
    \ifnum\@minsofar = \z@
      \@minsofar\number\@cite@temp
      \edef\@scan@copy{\@scan@cites}\else
    \ifnum\@cite@temp < \@minsofar
      \@minsofar\number\@cite@temp
      \edef\@scan@copy{\@scan@cites}\fi\fi\fi}\@tempcnta\@min
  \ifnum\@minsofar > \z@ 
    \advance\@tempcnta\@ne
    \@min\@minsofar
    \ifnum\@tempcnta=\@minsofar 
      \ifx\@h@ld\relax
        \edef \@h@ld{\@citea\csname b@\@scan@copy\endcsname}%
    \else \edef\@h@ld{\ifmmode{-}\else--\fi\csname b@\@scan@copy\endcsname}%
      \fi
    \else \@h@ld\@citea\csname b@\@scan@copy\endcsname
          \let\@h@ld\relax
  \fi 
\fi}%
\def\@citea{,\penalty\@highpenalty\,}}\@h@ld}{#1}}
\def\appendixname{Appendix}
\def\appendix{\par
  \def\pre@section{\appendixname{}}
  \setcounter{section}{1}
  \@addtoreset{equation}{section}
  \def\thesection{\Alph{section}}
  \def\theequation{\arabic{equation}}}

\makeatother
\begin{document}
\addtolength{\unitlength}{-0.5\unitlength}
\def\t{\theta}
\def\ov{\overline}
\def\a{\alpha}
\def\b{\beta}
\def\d{\delta}
\def\g{\gamma}
\def\wt{\widetilde}
\def\m{\overline m}
\def\f{\overline f}
\def\p{\overline p}
\def\n{\overline n}
\def\q{\overline q}
\def\w{\omega}
\def\ds{\displaystyle}
\def\s{\sigma}
\phantom{aa}
\hfill August, 1995

\vspace*{4cm}
\centerline{\bf $\Psi$ - VECTORS FOR THREE DIMENSIONAL MODELS.}
\vspace{1cm}
\centerline{
S.M. Sergeev\footnote{Branch Inst. for Nucl. Phys.,E-mail:
sergeev\_ms@mx.ihep.su},
G.E. Boos\footnote{E-mail: boos@mx.ihep.su},
V.V. Mangazeev\footnote{Australian National University,
E-mail: vvm105@phys.anu.edu.au}
and Yu.G. Stroganov\footnote{E-mail: stroganov@mx.ihep.su}}

\vspace{1cm}
\centerline{Institute for High Energy Physics,}
\centerline{Protvino, Moscow Region, Russia}

\vspace{1cm}
\centerline{Abstract}
In this paper we apply the method of $\psi$ - vectors to the
three dimensional statistical models. This method gives the
correspondence between the Bazhanov -- Baxter model and
its vertex version. Considering $\psi$ - vectors for the Planar model
we obtain its self -- duality.

\newpage

\section{Introduction}

Recently the vertex version of the Bazhanov -- Baxter model
(BBM) \cite{bb}
has been obtained in Ref. \cite{mss3}.
There we showed
the thermodynamic equivalence between the vertex model
and Interaction Round a Cube (IRC for the shortness) type BBM.
As the vertex model, it contains the Hietarinta's and Korepanov's
models \cite{hietarinta,korepanov}.

There are at least to types of vertex -- IRC correspondence.
The first one is the Wu - Kadanoff duality, which is valid
when it is possible  to regard the vertex variables of a weight
function as combinations of the IRC spins (see Ref. \cite{hietarinta}
for details). The second one is a three dimensional modification
of the $\psi$ vectors method \cite{abf,jap,cpm}.

In section 2 we fix the notations for the Tetrahedron equations of both
vertex and IRC types and give some evident notations for the $\psi$ -
vectors.
Necessary notations and definitions are given in section 3.
Next, in section 4 we give the explicit
expressions for the $\psi$ - vectors in BBM and
show the explicit equivalence between vertex and IRC types.
In section 5
we consider the Planar model \cite{mss2}
and find the $\psi$ - vectors for it. For this model
the vertex and the IRF type Boltzmann weights are the same.

Note that this paper does not contain explicit description of our
calculations, and for the details of the calculations see
the collection of $\omega$- hypergeometrical formulas
in the Appendix of Ref. \cite{mss3}.

\section{Tetrahedron equation and definition of $\psi$ - vectors.}

We consider two forms of the Tetrahedron equation (TE): vertex and IRC.
The vertex form is the following:
\begin{eqnarray}
&{\ds\sum_{k_1,k_2,k_3,\atop k_4,k_5,k_6}}
R^{k_1,k_2,k_3}_{i_1,i_2,i_3}
R'^{j_1k_4k_5}_{\phantom{,}k_1i_4i_5}
R''^{j_2j_4k_6}_{\phantom{,,}k_2k_4i_6}
R'''^{j_3j_5j_6}_{\phantom{,,,}k_3k_5k_6}
=&\nonumber\\
&={\ds\sum_{k_1,k_2,k_3,\atop k_4,k_5,k_6}}
R'''^{k_3,k_5,k_6}_{\phantom{,,,}i_3,i_5,i_6}
R''^{k_2k_4j_6}_{\phantom{,,}i_2i_4k_6}
R'^{k_1j_4j_5}_{\phantom{,}i_1k_4k_5}
R^{j_1j_2j_3}_{k_1k_2k_3}
&.\label{vte}
\end{eqnarray}
The IRC type TE is
\begin{eqnarray}
&\ds\sum_{d}
W(a_1|c_{12},c_{13},c_{14}|b_2,b_3,b_4|d)
W'(c_{12}|a_2,b_4,b_3|d,c_{24},c_{23}|b_1)&\nonumber\\
&\ds W''(b_4|c_{23},c_{13},d|b_2,b_1,a_3|c_{34})
W'''(d|b_1,b_2,b_3|c_{14},c_{24},c_{34}|a_4)=&\nonumber\\
&\ds = \sum_d
W'''(b_4|c_{23},c_{13},c_{12}|a_1,a_2,a_3|d)
W''(c_{12}\a_2,a_1,b_3|c_{14},c_{24},d|a_4)&\nonumber\\
&\ds W'(a_1|d,c_{13},c_{14}|b_2,a_4,a_3|c_{34})
W(d|a_2,a_3,a_4|c_{34},c_{24},c_{23}|b_1)&\label{wte}
\end{eqnarray}
This equation differs from the standard equation (2.2) in Ref. \cite{b}
by a simple reordering of the spins.

With each $R$ and $W$ we associate an oriented trihedron, described by
three ordered
dihedral angles $\t_1,\t_2,\t_3$, or, equivalently, by the corresponding
planar angles $a_1,a_2,a_3$. These trihedrons will be regarded
as the spectral arguments of the corresponding $R$ and $W$.
With TE we associate an (oriented) tetrahedron.
Complete solutions of TE are parameterized by six angles of the
tetrahedron (five of them are independent):
\begin{eqnarray}
&\ds (R,W) = (R,W)(\t_1,\t_2,\t_3),&\nonumber\\
&\ds (R',W') = (R,W)(\t_1,\t_4,\t_5),&\nonumber\\
&\ds (R'',W'') = (R,W)(\pi-\t_2,\t_4,\t_6)&\nonumber\\
&\ds (R''',W''') = (R,W)(\t_3,\pi-\t_5,\t_6).&        \label{angles}
\end{eqnarray}
The ordering of the dihedral angles is natural with respect to numbering
of the spaces and differs from that in the
standard equation (2.2) in Ref. \cite{b}.

For each vertex in (\ref{angles}) let $a^\sharp_i$ be the corresponding
planar angles:
\begin{eqnarray}
&\ds (\t_1,\t_2,\t_3)\rightarrow (a_1,a_2,a_3),&\nonumber\\
&\ds (\t_1,\t_4,\t_5)\rightarrow (a'_1,a'_2,a'_3),&\nonumber\\
&\ds (\pi-\t_2,\t_4,\t_6)\rightarrow (a''_1,a''_2,a''_3),&\nonumber\\
&\ds (\t_3,\pi-\t_5,\t_6)\rightarrow (a'''_1,a'''_2,a'''_3),&\label{linear}
\end{eqnarray}
and this correspondence will be implied below.

Vectors $\psi$ and $\ov\psi$ are defined as solutions of the following
equations:
\begin{eqnarray}\label{psi_eq}
&\ds \sum_{k_1,k_2,k_3} R^{k_1,k_2,k_3}_{i_1,i_2,i_3}
\psi_1(k_1|e,h,c,d)\psi_2(k_2|d,b,h,f)\psi_3(k_3|h,g,c,b)=&\nonumber\\
&\ds =\sum_{a}\psi_1(i_1|a,b,g,f)\psi_2(i_2|e,g,c,a)\psi_3(i_3|d,a,e,f)
W(a|e,f,g|b,c,d|h),&\nonumber\\
\end{eqnarray}

\begin{eqnarray}\label{ovpsi_eq}
&\ds \sum_{h} W(a|e,f,g|b,c,d|h)
\ov\psi_1(j_1|e,h,c,d)\ov\psi_2(j_2|d,b,h,f)
\ov\psi_3(j_3|h,g,c,b)=&\nonumber\\
&\ds =\sum_{k_1,k_2,k_3}
\ov\psi_1(k_1|a,b,g,f)\ov\psi_2(k_2|e,g,c,a)
\ov\psi_3(k_3|d,a,e,f) R^{j_1,j_2,j_3}_{k_1,k_2,k_3},&\nonumber\\
\end{eqnarray}
where $R$ and $W$ must obey  TEs.
The geometry demands the tetrahedron -- like parameterization
for $\psi$ and $\ov\psi$:
\begin{eqnarray}
&\ds (R,W) = (R,W)(\t_1,\t_2,\t_3),&\nonumber\\
&\ds (\psi_1,\ov\psi_1) = (\psi,\ov\psi)(\t_1,\t_4,\t_5),&\nonumber\\
&\ds (\psi_2,\ov\psi_2) = (\psi,\ov\psi)(\pi-\t_2,\t_4,\t_6)&\nonumber\\
&\ds (\psi_3,\ov\psi_3) = (\psi,\ov\psi)(\t_3,\pi-\t_5,\t_6).&
\label{psi_angles}
\end{eqnarray}
TEs and the equations for the psi vectors have two popular limits:
{\bf the static limit}, when in each $(R,W)(\t_1,\t_2,\t_3)$
$\t_1+\t_2+\t_3=\pi$, and {\bf the planar limit}, when
in each $(R,W)(a_1,a_2,a_3)$ $a_2=a_1+a_3$.

\section{Notations and definitions}

We tried to make the list of definitions and notations
shortest. Here we give only definition of $w$ function.
Its properties the Reader can find in Ref. \cite{mss3}.
Let
\begin{equation}
\w^{1/2}=\exp(\pi i/N).                                \label{2.1}
\end{equation}
Taking $p$ to be a point on a Fermat curve
$\Upsilon$, so that there defined three
complex numbers
$x(p),y(p),z(p)$, constrained by the
Fermat equation
\begin{equation}
x(p)^N+y(p)^N=z(p)^N,                                           \label{2.2}
\end{equation}
and $a$ to be an element of $Z_N$, define
\begin{equation}
\ds{w(p|a)\over w(p|0)}=\prod_{s=1}^{a}{y(p)\over z(p)-x(p)\w^s}. \label{2.3}
\end{equation}
The absolute value of $w(p|0)$ we define through
\begin{equation}
\prod_{a=0}^{N-1}w(p|a) = 1.         \label{2.4}
\end{equation}
Branches of $y(p)$ and $w(p|0)$ are arbitrary in general, but it is
convenient to choose them appropriately. Below all points $p$-s will be
defined so that (when it is possible, we will omit the argument $p$
for the shortness)
\begin{equation}\label{phreg}
\ds -2\pi/N<\mbox{Arg}(x/z)<0 \;\mbox{and}\; -\pi/N<\mbox{Arg}(y/z)<\pi/N.
\end{equation}
This subregion in $\Upsilon$ we call $\Upsilon_0$.
For $p\in\Upsilon_0$ we define $w(p|0)$ as follows
\begin{equation}\label{2.6}
\ds w(p|0) = \left({y\over z}\right)^{N-1\over 2}{1\over d(\w x/z)}=
\left({x\over y}\right)^{N-1\over 2}\Phi_0^{-1}d(z/x),
\end{equation}
where
\begin{equation}
\ds\Phi_0=\exp({i\pi (N-1)(N-2)\over 6N}).                      \label{2.7}
\end{equation}
and
\begin{equation}\label{2.8}
\ds d(x) = \exp\sum_{a=1}^{N-1} {a\over N}\log(1-x\w^a).
\end{equation}
It is implied in  (\ref{2.8}) that
\begin{equation}\label{2.9}
-\pi < \mbox{Im}\log(z) \leq \pi.
\end{equation}
Defined $w$ functions have the following property:
\begin{equation}
w(p|a)w(Op|-a)\wt\Phi(a)=1,\quad a\in Z_N,                    \label{2.10}
\end{equation}
where automorphism $O:\Upsilon_0\rightarrow\Upsilon_0$ is defined as
\begin{equation}
x(Op)=z(p),\; y(Op) = \w^{1/2}y(p),\; z(Op) = \w x(p),        \label{2.11}
\end{equation}
and
\begin{equation}
\wt\Phi(a)=\w^{a(a-N)/2}\exp({i\pi (N^2-1)\over 6N}),\;\;\;
\prod_{a\in Z_N}\wt\Phi(a)=1.                                 \label{2.12}
\end{equation}
There are a lot of several identities for the $w$ functions,
closely connected with the basic $q$ hypergeometric series.
Most useful of them
are the so called Star -- Square relation and $(\tau\rho)^2$
transformation. These identities are cumbersome enough, so we
do not write them here, and the Reader can find them
in Ref. \cite{mss3}.

\section{$\psi$ - vectors for BBM}

Recall the definition of the  vertex and IRC weights for  BBM.
The forms of the weights are taken from Refs. \cite{bb,mss3,kms}.
First, for given trihedron define four points $p_i=p_i(a_1,a_2,a_3)$:
\begin{eqnarray}
&\ds  x_{p_1}=\w^{-1/2}\exp (i{a_3\over N})\sqrt[N]{\sin\b_1\over\sin\b_2},\;
y_{p_1} = \exp ( i{\b_1\over N})\sqrt[N]{\sin a_3\over\sin\b_2};&\nonumber\\
&\ds  x_{p_2}=\w^{-1/2}\exp (i{a_3\over N})\sqrt[N]{\sin\b_2\over\sin\b_1},\;
y_{p_2} = \exp ( i{\b_2\over N})\sqrt[N]{\sin a_3\over\sin\b_1};&\nonumber\\
&\ds  x_{p_3}=\w^{-1}\exp (i{a_3\over N})\sqrt[N]{\sin\b_3\over\sin\b_0},\;
y_{p_3} = \exp ( -i{\b_3\over N})\sqrt[N]{\sin a_3\over\sin\b_0};&\nonumber\\
&\ds  x_{p_4}=\w^{-1}\exp (i{a_3\over N})\sqrt[N]{\sin\b_0\over\sin\b_3},\;
y_{p_4} = \exp ( -i{\b_0\over N})\sqrt[N]{\sin a_3\over\sin\b_3};&\nonumber\\
&z_{p_i}=1,\,i=1,2,3,4;&                              \label{p_points}
\end{eqnarray}
where the linear excesses are
\begin{equation}
\ds\b_0=\pi-{a_1+a_2+a_3\over 2},\;\b_i={a_j+a_k-a_i\over 2}.
                                                       \label{excesses}
\end{equation}
Let $\rho_k$ be the normalization factors:
\begin{equation}
\ds \rho_k = {1\over N}\left({\sin a_k\over 2\cos\b_0/2\dots\cos\b_3/2}
\right)^{N-1\over N}.
\end{equation}
The vertex weight is given by
\begin{equation}
\ds R^{j_1,j_2,j_3}_{i_1,i_2,i_3}=
\ds\d_{j_2+j_3,i_2+i_3}\w^{j_3(j_1-i_1)}
\rho_3{w(p_1|i_1-i_2)w(p_2|j_1-j_2)\over
 w(p_3|i_1-j_2)w(p_4|j_1-i_2)}.                       \label{v_weight}
\end{equation}
Introducing the normalization factor $\rho_3$, we make $R$ symmetrical
with respect to cube symmetry group (see Ref. \cite{mss3})
and restore Bazhanov -- Baxter's  normalization of the model
(see Ref. \cite{bb}).

Define other four points $q_i(a_1,a_2,a_3)=Op_i(a_1,a_3,a_2)$
The IRC weight is given by
\begin{eqnarray}
&\ds W(a|e,f,g|b,c,d|h|\t_1,\t_2,\t_3) =&\nonumber\\
&\ds =\rho_2\sum_\s{w(q_4|f-a+\s)w(q_3|h-c+\s)\over
w(q_1|d-e+\s) w(q_2|b-g+\s)}\w^{\ds\s(e+g-a-c)}.&\label{w_weight}
\end{eqnarray}
This weight differs from the Bazhanov -- Baxter weight,
the correspondence is
\begin{equation}
\ds W_{B}(a|e,f,g|b,c,d|h|\t^B_1,\t^B_2,\t^B_3) =
W(a|f,g,e|c,d,b|h|\t_1,\t_2,\t_3)\label{corresp_w}
\end{equation}
where
\begin{equation}
\t_1=\t^B_2,\;\;\t_2=\t^B_3,\;\;\t_3=\t^B_1.\label{corresp_t}
\end{equation}

The formulae for $\psi$ and $\ov\psi$ are:
\begin{eqnarray}
&\ds \psi(\s|e,h,c,d) = {w(s|\s+e-c)\over w(t|\s+d-h)}
\w^{\s(h-c)},&\nonumber\\
&\ds \ov\psi (\s|a,b,g,f) = {w(s'|\s+f-b)\over w(t'|\s+a-g)}
\w^{\s(g-b)},\label{psi}
\end{eqnarray}
where if these $(\psi,\ov\psi)=(\psi,\ov\psi)(\t_1,\t_2,\t_3)$ then
\begin{eqnarray}
&\ds s = q_4(a_2,\pi-a_3,\pi-a_1),\;\;
t = q_1(a_2,\pi-a_3,\pi-a_1),&\nonumber\\
&\ds s' = q_3(a_2,\pi-a_3,\pi-a_1),\;\;
t' = q_2(a_2,\pi-a_3,\pi-a_1),&\label{psi_param}
\end{eqnarray}

Note that a natural IRF -- type $L$ -- operator
\begin{equation}
\ds L(a|e,f,g|b,c,d|h)_{\t_1,\t_2,\t_3}=
\rho_1\sum_\s\psi(\s|e,h,c,d)\ov\psi(\s|a,b,g,f)\label{L1_op}
\end{equation}
is equivalent to the $(\tau\rho)^2$ transformed weight $W$
(\ref{w_weight}).
The analogous vertex type $L$ -- operator
\begin{equation}
\ds L_{i,c-e,e-d}^{j,h-d,c-h} = \rho_1\psi(i|e,h,c,d)\ov\psi(j|e,h,c,d)
\end{equation}
is also equivalent to the appropriately transformed vertex weight
$R$ (\ref{v_weight}).

The proof of equations (\ref{psi_eq},\ref{ovpsi_eq}) are simple.
For the left hand side of (\ref{psi_eq}) one has to make $(\tau\rho)^2$
over $k_1$ and after this the Star -- Square summation formula over $k_2$
(see Ref. \cite{mss3} for the meaning of these charms).
In the right hand side one has to sum over the spin $a$. The final
expressions are to coincide.

\section{$\psi$ - vectors for the Planar Model.}

The Planar model is considered in Ref. \cite{mss2}. Recall the
definition of it.

Define another four points for the trihedron $(a_1,a_2,a_3)$:
\begin{equation}
\ds r_i=\left(\exp(-i{\b_i\over N}),\w^{1/4}\sqrt[N]{2\sin\b_i},
\exp(i{\b_i\over N})\right).
\end{equation}
Consider the planar limit of the TE, i. e. the case when $\b_2=0$ for each
weight. The vertex weight is
\begin{equation}
\ds R^{j_1,j_2,j_3}_{i_1,i_2,i_3}=
\ds\d_{j_2,i_1+i_3}\d_{i_2,j_1+j_3}\w^{j_1(i_3-j_3)}
{w(r_1|i_3-j_3)w(r_3|i_1-j_1)\over
 w(Or_0|j_2-i_2)}.                       \label{v_plm}
\end{equation}
IRC form of this weight is
\begin{eqnarray}
&\ds W(a|e,f,g|b,c,d|h)=&\nonumber\\
&\ds\w^{(h-e)(a-d-g+h)}
{w(r_1|a-d-g+h)w(r_3|b-a-h+e)\over
 w(Or_0|b-d-g+e)}.&                       \label{w_plm}
\end{eqnarray}
The Planar model is self -- dual in the sense that
\begin{equation}
\ds R^{h-e,b-d,g-h}_{b-a,g-e,a-d}=
W(a|e,f,g|b,c,d|h).
\end{equation}

$\psi$ - vectors for the weights (\ref{v_plm}) and (\ref{w_plm}) are
\begin{eqnarray}
&\ds \psi(\s|a,b,c,d) = w(v|\s+a-b)\w^{\s(d-b)}\phi(a,b,c,d),&\nonumber\\
&\ds \ov\psi (\s|a,b,c,d) = {\w^{\s(c-a)}\over w(u|\s+a-b)}\ov\phi(a,b,c,d),&
\label{psi_plm}
\end{eqnarray}
where for each $\phi$ and $\ov\phi$ there are two independent choices:
\begin{eqnarray}
&\phi(a,b,c,d)=\w^{(a-b)(a-c)}\;\mbox{or}\;\w^{(a-b)(d-b)},&\nonumber\\
&\ov\phi(a,b,c,d)=\w^{(a-b)(b-d)}\;\mbox{or}\;\w^{(a-b)(c-a)}.&\label{psi_ph}
\end{eqnarray}
It is useful to choose the second expressions for the phases $\phi$ and
$\ov\phi$.
The equations for the $\psi$ - vectors are equivalent to a
pair of $(\tau\rho)^2$ transformations.
According to the scheme (\ref{psi_angles}) the arguments of $\psi$
 vectors are (omitting the middle argument $a_2=a_1+a_3$):
\begin{eqnarray}\label{vu}
&\ds v(a_1,a_3)=\left(
\exp(-i{a_1\over N})\sqrt[N]{\sin a_3\over\sin a_2},
\exp(i{a_3\over N})\sqrt[N]{\sin a_1\over\sin a_2},1\right),&
\nonumber\\
&\ds u(a_1,a_3)=\left(
\w^{-1}\exp(i{a_1\over N})\sqrt[N]{\sin a_3\over\sin a_2},
\exp(-i{a_3\over N})\sqrt[N]{\sin a_1\over\sin a_2},1\right).&
\end{eqnarray}
The $\psi$ - vectors are in the same time $L$
operators for the $R$ matrix (\ref{v_plm}):
\begin{equation}
\ds L^{j_1,j_2,j_3}_{i_1,i_2,i_3} =
\d_{j_1,j_3-i_2}\d_{j_1,i_3-j_2} w(v|i_1-j_1)\w^{j_2(i_1-j_1)}
\end{equation}
and analogically for $\ov\psi$.
With the help of psi vectors one can restore the $W$ weight (\ref{w_plm}).
\begin{eqnarray}\label{dec_psi}
&\ds n_1\sum_\s\psi_v(\s|e,h,c,d)\ov\psi_u(\s|a,b,g,f)=&\nonumber\\
&\ds W(a|e,f,g|b,c,d|h)\w^{-(a-b)(d-h)-(a-g)(h-e)}.&
\end{eqnarray}
The phase factor in the right hand side is the gauge factor
for TE and the normalization factor
\begin{equation}
\ds n_k^{-1} = \sqrt{N}\exp ( i\pi{N^2-1\over 12N})
\left(2{\sin a_i\sin a_j\over\sin a_k}\right)^{N-1\over 2N}
\end{equation}

\section{Discussion}

In two dimensional statistical systems the application of the $\psi$
vectors method gives the excellent results. First, the IRF type models,
corresponding to the $R$ matrices of the simple Lie algebras,
admit the elliptic deformation even if the $R$ matrices are
trigonometric \cite{jap}.  Second, the $\psi$ vectors for
the cyclic representations of ${\cal U}_q(A_n)$ allow one
to construct new models \cite{cpm}. Contrary to this situation
in $D=2$, we still have not succeeded in any
application of the three dimensional $\psi$ vectors.

In three dimensions we have three types of Boltzmann weights.
Two of them, vertex and IRC, were considered above. The third
type Boltzmann weights with twelve spin variables are connected with
the scattering straight strings \cite{z}. Taking into account
this third type of the spin structure, one can modify the given scheme
of the $\psi$ vectors in different ways. Perhaps, in the case
when $\psi$ vectors for some modified scheme exist, they would
give us a method of constructing new models.

\noindent
{\bf Acknowledgements}

\noindent
We would like to thank J.- M. Maillard, J. Hietarinta and I. Korepanov
for useful discussions. One of us (S.S.) is also grateful to the
organizing committee of the International Conference of Mathematical
Physics (Chelyabinsk, july 1995) for hospitality.
This research has been partially supported by National Science
Foundation Grant PHY -- 93 -- 07 -- 816, by International Science
Foundation  Grant RMM300, by Russian Foundation of Fundamental
Research Grant 95 -- 01 -- 00249.

\end{document}